\input{amssym.def}
\input{amssym}

\newtheorem{cor}{Corollary}

\newtheorem{lemma}{Lemma}

\newtheorem{thm}{Theorem}

\newtheorem{claim}{Claim}

\newcommand{\be}{\begin{enumerate}}
\newcommand{\ee}{\end{enumerate}}
\newcommand{\bq}{\begin{quote}}
\newcommand{\eq}{\end{quote}}
\newcommand{\forces}{\Vdash}

\newcommand{\wh}{\widehat{}\ }

\documentstyle[12pt
]{article}

\begin{document}
\author{Boban Velickovic\\ and\\ W. Hugh Woodin}
\title{Complexity of the reals of inner models of set theory}
\date{}
\maketitle

\section*{Introduction}

The usual definition of the set of constructible reals ${\Bbb R}^L$ is $\Sigma ^1_2$. 
This set can have a simpler definition if, for example, it is countable or if
every real is in $L$. 
Martin and Solovay [MS1] have shown that if MA$_{\aleph _1}$ holds and
there is a real $r$ such that $\aleph _1^{L[r]}=\aleph _1$ 
then every set of reals of size $\aleph _1$ is  co-analytic.
Thus by a ccc forcing over a universe of $V=L$ we can obtain a 
universe of set theory
in which ${\Bbb R}^L$ is an uncountable co-analytic set yet not every real is in $L$. 
The results of this paper were motivated by a question of H. Friedman [Fr, problem 86]
who asked if ${\Bbb R}^L$ can be analytic or even Borel
in a nontrivial way, that is both uncountable and not equal to the set
of all reals. 
There is a companion question due to K. Prikry whether ${\Bbb R}^L$
could contain a perfect set and not be equal to the set of all reals.
Clearly a positive answer to the first question 
would also imply a positive answer to the second one.

The main result of this paper is a negative answer to Friedman's question.
In fact we prove that if $M$ is an inner model of set theory and the set ${\Bbb R}^M$
of reals in $M$ is analytic then either all reals are in $M$
or else $\aleph _1^M$ is countable. 
Since the cardinality of ${\Bbb R}^L$ is $\aleph _1^L$
this implies the desired result in the case $M=L$. 
We also show that in the context of large cardinals this result can be extended to 
projective sets in place of analytic sets.
However, the conclusion of the main theorem cannot be strengthened
to say that either all reals are in $M$ or else the continuum of $M$
is countable. We produce a pair of generic extensions $W$ and $V$  of $L$
such that $W\subseteq V$, the reals of $W$ form an uncountable $F_{\sigma}$
set in $V$, and yet not all reals from $V$ are in $W$.

In relation to Prikry's problem we show that if an inner model $M$
contains a superperfect set of reals then it contains all reals. 
The proof is based on a construction of a recursive coloring 
of triples of reals into $2^{\omega}$ such that for any superperfect
set $P$  the triples from $P$ obtain all colors.
A similar partition was used by Gitik [Gi] who showed that if $V$ 
is a universe of set theory and $r$ is a real not in $V$ then 
the set of countable subsets of $\omega _2$ in $V[r]$ which are 
not in $V$ form a stationary set in $[\omega _2]^{\aleph _0}$.
 It was observed by the first
author in [Ve] that this implies that if the Semi Proper Forcing Axiom (SPFA)
holds and $M$ is an inner model of set theory such that $\aleph _2^M=\aleph _2$
then all reals are in $M$. 
In the positive direction of Prikry's problem we give an example 
of two generic extensions $V$ and $W$ of $L$ such that $W\subseteq V$,
$\aleph _1^W=\aleph _1^V$, and there is a perfect set in $V$ consisting 
of reals from $W$. 

The paper is organized as follows. In \S 1 we present the coloring of
triples of reals described above. It uses oscillations of reals numbers,
a technique commonly used in the construction of examples to negative partition relations 
(see for example [To]). We then use it to prove a special case of
the main theorem in the case of the constructible universe. 
Although this proof, which uses Jensen's Covering Lemma, is subsumed
by Theorem 3 we present it since it may have some interest
of its own.

In \S2  we prove a kind of regularity property for $\bf \Sigma$$_1^1$ sets
saying that if an analytic set $A$ contains codes for all countable ordinals
then every real is hyperarithmetic in a finite sequence of elements of
$A$. From this our main result follows easily.  We then extend this
to higher levels of the projective hierarchy under appropriate 
large cardinal assumptions or projective determinacy.

Section \S 3 contains examples of pairs of models of set theory which
show that the above results are in some sense best possible.
We prove that it is possible to have an inner model of set theory $W$ 
whose reals form an uncountable $F_{\sigma}$ set and yet not all reals
belong to $W$. Then necessarily $\aleph _1^W$ is countable. 
However it is possible to have $\aleph _1^W=\aleph _1$ if we only
require that $W$ contains a perfect set of reals. 

Finally in  \S 4 we prove assuming AD + $V=L({\Bbb R})$ that   if $M$ 
is an inner model of ZF containing a Souslin prewellordering of reals
of length $\aleph _1^V$ then all reals are in $M$. 
This result has some consequences in the theory of cardinals 
in $L({\Bbb R})$ under the axiom of determinacy. 
Some assumption on the prewellordering in the above result is necessary.
We prove in ZF alone that assuming there is a nonconstructible real
 there is an inner model $M$ of ZF containing a prewellordering of
reals of length $\aleph _1^V$ and such that not all reals belong to $M$. 

Our notation is fairly standard or self-explanatory. For all undefined
notions see [Ku]. 
For an index set $I$ we shall let ${\cal C}(I)$ denote the usual forcing
for adding $I$ Cohen reals. Thus conditions in ${\cal C}(I)$ are
finite partial functions from $\omega \times I$ to $\{ 0,1\}$ 
and the order is inclusion.

\section{Coloring triples of reals}

We now present the coloring of triples of reals described in the introduction.
First we make some relevant definitions. 
We  identify the set of reals ${\Bbb R}$ with the set 
$(\omega)^{\omega}$ of all infinite increasing sequences of natural numbers.
We shall let $\leq _*$ denote the ordering of eventual dominance on
$(\omega)^{\omega}$.
We also let $(\omega)^{<\omega}$ denote the set of all finite increasing
sequences of natural numbers. Then $(\omega)^{<\omega}$ forms a tree 
under inclusion.  
Given a subtree $T$ of $(\omega)^{<\omega}$ we say that a node $s\in T$ is 
{\em $\omega$-splitting} if the set $\{ k: s\ \wh  k\in T\}$ is infinite. 
$T$ is called {\em superperfect}  if above every node
$s\in T$ there is a node  $t\in T$ which is $\omega$-splitting.
A subset $P$ of $(\omega)^{\omega}$ is called {\em superperfect} 
if the set $T$ of all finite initial  segments of members of $P$ forms 
a superperfect tree.

\begin{thm}{There is a partial recursive  function 
$o:{\Bbb R}^3\rightarrow \{ 0,1\} ^{\omega}$ such that for 
every superperfect set $P$ $o '' (P^3) =\{ 0,1\}^{\omega}$.}
\end{thm}

\noindent PROOF: Given distinct reals $x,y,z \in (\omega)^{\omega}$ let 
$$
O(x,y,z)=\{ n: z(n-1)\leq x(n-1),y(n-1)\ \mbox{and}\ x(n),y(n)<z(n) \}.
$$
$o(x,y,z)$ will be defined iff $O(x,y,z)$ is infinite.
If $O(x,y,z)$ is infinite let $\{n_k:k<\omega\}$ be the increasing enumeration 
of its members.
Define $o(x,y,z)$ to be the function $\alpha :\omega\rightarrow \{ 0,1\}$ 
where for every $k < \omega$,
$$
\alpha (k)=0\ \mbox{iff}\ x(n_k)\leq y(n_k).
$$
We show that if $P$ is a superperfect subset of $(\omega)^{\omega}$ and 
$\alpha \in \{ 0,1\}^{\omega}$ there are $x,y,z \in P$ such that 
$o(x,y,z)=\alpha$. 

Thus, fix such a superperfect set $P$ and let $T$ be the tree of all 
initial segments of elements of $P$. 
We define inductively strictly increasing sequences of $\omega$-splitting nodes
of $T$ $\{ x_k:k<\omega \}$, $\{y_k:k<\omega \}$, $\{ z_k:k<\omega \}$ as 
follows. 
The lengths of $x_k,y_k,$ and $z_k$ will be $l_k,m_k,$ and $n_k$ respectively.
We will have $n_k<l_k,m_k$ and $z_k(n_k-1)<x_k(n_k-1),y_k(n_k-1)$. 
At stage $k$ we use the fact that $z_k$ is an $\omega$-splitting node
of $T$ to find an integer
$i>x_k(l_k-1),y_k(m_k-1)$ such that $z_k\ \wh i \in T$. 
We then let $z_{k+1}$ be any $\omega$-splitting node of $T$ 
of length $n_{k+1}>l_k,m_k$
extending $z_k\ \wh i$. 
Now look at $\alpha (k+1)$.
Let as assume for definiteness that it is equal $0$.
Since $x_k$ is a splitting node of $T$ we can find an integer 
$j>z_{k+1}(n_{k+1}-1)$ 
such that $x_k\ \wh j\in T$. Now let $x_{k+1}$ be any $\omega$-splitting
node of $T$ extending $x_k\ \wh j$ of length $l_{k+1}>n_{k+1}$.
Finally find some integer $h>x_{k+1}(l_{k+1}-1)$ such that $y_k\ \wh h\in T$ and
let $y_{k+1}$ be any $\omega$-splitting node of $T$ extending 
$y_k\ \wh h$ of length $m_{k+1}\geq l_{k+1}$. 
If $\alpha (k+1)=1$ then reverse the construction of $x_{k+1}$ and $y_{k+1}$.
At the end we let 
$x=\bigcup \{x_k: k<\omega \}$, $y=\bigcup \{ y_k:k<\omega \}$, and
$z=\bigcup \{ z_k:k<\omega \}$. We then have $o(x,y,z)=\alpha$.
\hspace{.25in} $\Box$

\medskip

\begin{cor}{ Let $V$ and $W$ be models of set theory such that $W$ is a 
subuniverse of $V$. If $V$ contains a superperfect tree $T$ all of whose
branches lie in $W$ then $V$ and $W$ have the same reals. \hspace{.25in} $\Box$}
\end{cor}

\medskip

Given a limit ordinal $\lambda$ let $(\lambda)^{<\omega}$ denote the set of
finite increasing sequences of ordinals $<\lambda$ and let $(\lambda)^{\omega}$ 
denote the set of all increasing $\omega$-sequences in $\lambda$. 
Say that a subtree $T$ of $(\lambda)^{<\omega}$ is $\lambda$-{\em superperfect} if
for every node $s\in T$ there is $t\in T$ extending $s$ such that
the set $\{ \alpha : t\ \wh \alpha \in T\}$ is cofinal in $\lambda$. 
Say that a subset $P$ of  $(\lambda)^{\omega}$ is $\lambda$-{\em superperfect}
if  the set of finite initial segments of members of $P$ forms a a $\lambda$-superperfect
tree. The same construction as in Theorem 1 gives a coloring 
$o_{\lambda}:((\lambda)^{\omega})^3\rightarrow \{ 0,1\}^{\omega}$ such
that for any $\lambda$-superperfect set $P\subseteq (\lambda)^{\omega}$ 
$o_{\lambda}''P^3=\{ 0,1\}^{\omega}$. 
Moreover if $x,y,z\in (\lambda)^{\omega}$ then $o_{\lambda}(x,y,z)\in L[x,y,z]$.

\medskip

\begin{thm}{Suppose ${\Bbb R}^L$ is an uncountable analytic set.
Then every real is constructible.}
\end{thm}

\medskip

\noindent PROOF:
By an old result of Hurewitz  every analytic subset of $(\omega)^{\omega}$
is either bounded under $\leq _*$ or contains a superperfect subset. 
Thus in order to complete the proof of Theorem 2 it suffices to
establish that under the assumptions of the theorem 
${\Bbb R}^L$ is unbounded in $(\omega)^{\omega}$ under $\leq _*$.
Note that in this case $\aleph _1^L=\aleph _1$ and hence 
$0^{\#}$ does not exist. 
Thus we can make use of Jensen's Covering Theorem [DJ].

Let us fix a subtree $T$ of $(\omega \times \omega)^{<\omega}$ such
that ${\Bbb R}^L=p[T]$. For any model $M$ of set theory containing
$T$ we denote the projection of $T$ in $M$ by $A^M$. 
Thus we have ${\Bbb R}^L=A^V$.
Fix a sufficiently large regular cardinal $\theta$ and let 
$N$ be a countable elementary submodel  of $H_{\theta}$, the collections
of sets hereditarily of size $<\theta$, which contains $T$ and let 
$M$ be the transitive collapse of $N$.  Then, by elementarity, we
have that ${\Bbb R}^{L^M}=A^M$. 

Work for a moment in $M$. 
Fix a singular cardinal $\kappa$. Then by the Covering Theorem applied
in $M$ $\kappa$ is also singular in $L^M$. 
Let $\lambda =\mbox{cof}^{L^M}(\kappa)$ and fix a sequence 
$\langle \kappa _{\xi}:\xi <\lambda \rangle \in L^M$ of $L^M$-regular 
cardinals converging to $\kappa$ such that 
$\kappa _0 >(\lambda ^{++})^M$.  We shall consider the orderings 
$\leq$ and $\leq ^*$ on $\Pi _{\xi <\lambda}\kappa _{\xi}$ 
of everywhere dominance and eventual dominance respectively.
By a straightforward application of the Covering Theorem we have
that $\Pi _{\xi <\lambda}\kappa _{\xi}\cap L^M$ is cofinal
in $\Pi _{\xi <\lambda}\kappa _{\xi}$ under $\leq$.  Also note that
$\Pi _{\xi <\lambda}\kappa _{\xi}$ is $<\kappa ^+$-directed under $\leq ^*$.
Let ${\cal P}$ denote $Coll(\aleph _0,\kappa)$, the standard poset for collapsing $\kappa$ to
$\aleph _0$ with finite conditions. Note that ${\cal P}$ has size $\kappa$.
Back in $V$ pick an $M$-generic filter $G$ over ${\cal P}$. 

\begin{claim}{$(\Pi _{\xi <\lambda}\kappa _{\xi})^M$ is unbounded in 
$(\Pi _{\xi <\lambda}\kappa _{\xi})^{M[G]}$ under $\leq ^*$.}
\end{claim}

\noindent PROOF: Suppose otherwise and work in $M$. 
Let $\tau$ be a ${\cal P}$-name
for a function in $\Pi _{\xi <\lambda}\kappa _{\xi}$ which 
eventually dominates all ground model functions.  For $p\in {\cal P}$ 
and $\mu <\lambda$ let 
$$
X_{p,\mu}=\{ f\in \Pi _{\xi <\lambda} \kappa _{\xi}: 
p\forces f(\eta)\leq \tau (\eta),\ \mbox{for all}\ \eta \geq \mu \}.
$$
Then the sets $X_{p,\mu}$, for $p \in {\cal P}$ and $\mu <\lambda$,
 cover $\Pi _{\xi <\lambda}\kappa _{\xi}$.
Let $g_{p,\mu}\in \Pi _{\xi <\lambda}\kappa _{\xi}$ be defined as 
follows. For $\eta <\mu$ let $g_{p,\mu}(\eta)=0$. For $\eta \geq \mu$ let
$$
g_{p,\mu}(\eta) =\mbox{sup}\{f(\eta): f\in X_{p,\mu}\}.
$$
Then each $g_{p,\mu}$ is well-defined and eventually dominates all 
$f \in X_{p,\mu}$. Since $\Pi _{\xi <\lambda}\kappa _{\xi}$ is
$<\kappa ^+$-directed under $\leq ^*$ 
it follows that there is a single function 
$g\in \Pi _{\xi <\lambda}\kappa _{\xi}$ which eventually dominates all
the $g_{p,\mu}$. Contradiction. \hspace{.25in} $\Box$

\medskip

Now in $M[G]$ find a linear ordering $<_{\lambda}$ such that 
$\langle \omega , <_{\lambda}\rangle$
 is isomorphic to $\langle \lambda ,<\rangle$ and such that
$<_{\lambda}$ coded in some reasonable way as an element of 
$(\omega)^{\omega}$ belongs to $A^{M[G]}$. 
Since $\lambda$ is countable in the true $L$ we can find such a 
linear order $<_{\lambda}$ in $A^V={\Bbb R}^L$ and hence by absoluteness
of $\bf \Sigma$$_1^1$ formulas between $V$ and $M[G]$ we can find 
such $<_{\lambda}$ in $A^{M[G]}$. 
Let now 
$e:\langle \omega ,<_{\lambda}\rangle \rightarrow \langle \lambda ,<\rangle$
be the unique isomorphism. By a similar argument we can find linear orderings
$<_n$ on $\omega$ such that $\langle \omega ,<_n\rangle$ is isomorphic to 
$\langle \kappa _{e(n)},<\rangle$ and the sequence 
$\langle <_n:n<\omega \rangle$
coded in some reasonable way is in $A^{M[G]}$. 
It follows that this sequence is in $L$ although of course not in 
$L^M$. Let $e_n:\langle \kappa _{e(n)},<\rangle \rightarrow 
\langle \omega ,<_n \rangle$ be the unique
isomorphism. Define a map 
$\varphi :\Pi _{\xi <\lambda}\kappa _{\xi} \rightarrow (\omega)^{\omega}$
as follows.
$$
\varphi (f)(n)= \Sigma _{i\leq n}e_i(f(e^{-1}(i))).
$$
Then a code of $\varphi$ exists in both $M[G]$ and $L$.
If $f\in \Pi _{\xi <\lambda}\kappa _{\xi}\cap {L^M}$ then 
$\varphi (f)$ is in $L$, hence in $A^V$, hence also in $A^{M[G]}$. Moreover 
$\varphi ''(\Pi _{\xi <\lambda}\kappa _{\xi}\cap L^M)$ is unbounded in 
${\Bbb R}^{M[G]}$ under $\leq _*$. Hence so is  $A^{M[G]}$.
Now by absoluteness of $\Sigma$$_2^1(T)$ formulas between $M$ and
$M[G]$ it follows that $A^M$ is unbounded in ${\Bbb R}^M$. 
Since $M$ is elementary equivalent to $H_{\theta}$ it follows
that $A^V$ is unbounded in ${\Bbb R}$. This finishes the proof
of Theorem 2. \hspace{.25in} $\Box$.

\section{Main Theorem}

In this section we prove the main result of this paper. We start with a lemma 
establishing a kind of regularity property for analytic sets of reals. 

\begin{lemma}{Suppose that $A$ is an analytic set such that 
sup$\{ \omega _1^{CK,x}: x\in A \}=\omega _1$. Then every real 
is hyperarithmetic in a quadruple of elements of $A$.}
\end{lemma}

\medskip

\noindent PROOF: Let $T \subset (\omega \times \omega)^{<\omega}$ 
be a tree such that $A=p[T]$. Note that the statement 
that sup$\{ \omega _1^{CK,x}:x\in p[T]\}$ is $\Pi _2^1(T)$ and
thus absolute. 

For an ordinal $\alpha$ let $Coll(\aleph _0,\alpha)$
be the usual collapse of $\alpha$ to $\aleph _0$ using finite
conditions. 
Let ${\cal P}$ denote $Coll(\aleph _0,\aleph _1)$.
 If $G$ is $V$-generic over ${\cal P}$, 
by Shoenfield's absoluteness theorem,
in $V[G]$ there is $x \in p[T]$ such that $\omega _1^{CK,x}>\omega _1^V$. 
In $V$ fix a name $\dot{x}$ for $x$ and a name $\sigma$ for a cofinal 
$\omega$-sequence in $\omega _1^V$ such that the maximal condition
in ${\cal P}$ forces that $\dot{x}\in p[T]$ and $\sigma \in L[\dot{x}]$. 

\medskip

\noindent {CLAIM 1:} For every $p \in {\cal P}$ there is $k<\omega$
such that for every $\alpha <\omega _1$ there is $q\leq p$ such that
$q\forces \sigma (k)>\alpha$. 

\medskip

\noindent PROOF: Assume otherwise and fix $p$ for which the claim is false.
Then for every $k$ there is $\alpha _k<\omega _1$ such that
$p \forces \sigma (k)<\alpha _k$. Let $\alpha =\sup \{ \alpha _k:k<\omega \}$.
Then $p\forces \mbox{ran}(\sigma)\subset \alpha$, contradicting the fact
$\sigma$ is forced to be cofinal in $\omega _1$. \hspace{.25in} $\Box$. 

\medskip

Let ${\cal Q}$ denote $Coll(\aleph _0,\aleph _2)$ as defined in $V$. 
 Suppose $H$ is $V$-generic over ${\cal Q}$. 
Work for a moment in $V[H]$. If $G$ is a $V$-generic filter over ${\cal P}$
let $\sigma _G$ denote the interpretation of $\sigma$ in $V[G]$.
Let $B$ be the set of all $\sigma _G$ where $G$ ranges over all $V$-generic
filters over ${\cal P}$. 

\medskip
\noindent {CLAIM 2:} $B$ contains an $\omega _1^V$-superperfect set
 in $(\omega _1^V)^{\omega}$. 

\medskip

\noindent PROOF: Let $\{ D_n:n<\omega \}$ be an enumeration of all dense subsets
of ${\cal P}$ which belong to the ground model.   For each
 $t\in (\omega _1^V)^{<\omega}$
we define a condition $p_t$ in the regular open algebra of ${\cal P}$
as computed in $V$ and $s_t\in (\omega _1^V)^{<\omega}$ inductively on the length
of $t$ such that

\begin{enumerate}
\item $p_t\in D_{lh(t)}$
\item $p_t \forces s_t \subset \sigma$
\item if $t\leq r$ then $p_r\leq p_t$ and $s_t \subset s_r$
\item if $t$ and $r$ are incomparable then $s_t$ and $s_r$ are incomparable
\item for every $t$ the set $\{ \alpha : \mbox{there is}\ q \leq p\ q\forces s_t\wh \alpha \subset \sigma\}$
is unbounded in $\omega _1^V$.
\end{enumerate}

Suppose $p_t$ and $s_t$ have been defined. Using 4. choose in $V$ 
a 1-1 order preserving function $f_t:\omega _1^V\rightarrow \omega _1^V$ 
and for every $\alpha$  $q_{t,\alpha}\leq p_t$ such that
 $q_{t,\alpha}\forces s_t\wh f_t({\alpha})\subset \sigma$.
By extending $q_{t,\alpha}$ if necessary we may assume that it belongs to $D_{lh(t)+1}$.
Now by applying Claim 1 we can find a condition $p\leq  q_{t,\alpha}$ 
and $k>lh(s_t)+1$ such that for some $s\in (\omega _1^V)^k$\
$p\forces s\subset \sigma$ and for every  $\gamma <\omega _1^V$
there is $q\leq p$ such that $q\forces \sigma (k)>\gamma$. 
Let then $s_{t\  \wh \alpha}=s$ and $p_{t\ \wh \alpha}=p$. 
This completes the inductive construction.

Now if $b\in (\omega _1^V)^{\omega}$ then $\{ p_{b\restriction n}:n<\omega\}$ 
generates a filter $G_b$ which is $V$-generic over ${\cal P}$.
The interpretation of $\sigma$ under $G_b$ is $s_b=\bigcup _{n<\omega}s_{b\restriction n}$.
Since the set $R=\{ s_b:b\in (\omega _1^V)^{\omega}\}$ is $\omega _1^V$-superperfect
this proves Claim 2. \hspace{.25in} $\Box$

Now using the remark following the proof of Theorem 1 for any real $r\in \{ 0,1\}^{\omega}$ 
we can find $b_1,b_2,b_3 \in (\omega _1^V)^{\omega}$
 such that $r \in L[s_{b_1},s_{b_2},s_{b_3}]$.
Let $x_i$ be the interpretation of $\dot{x}$ under $G_{b_i}$.
Then it follows that $x_i\in p[T]$ and $s_{b_i}\in L[x_i]$, for $i=1,2,3$. 
Thus $r\in L[x_1,x_2,x_3]$.
Pick a countable ordinal $\delta$ such that $r\in L_{\delta}[x_1,x_2,x_3]$.
Using the fact that in $V[H]$ 
sup$\{ \omega _1^{CK,x}:x\in p[T]\} =\aleph _1$ 
we can find $y\in p[T]$ such that $\omega _1^{CK,y}>\delta$. 
Then we have that $r$ is $\Delta _1^1(x_1,x_2,x_3,y)$.
Note that the statement that there are $x_1,x_2,x_3,y \in p[T]$ such that
$r\in \Delta _1^1(x_1,x_2,x_3,y)$ is $\Sigma _2^1(r,T)$. Thus for any $r \in V$,
by Shoenfield absoluteness again, it must be true in $V$. This proves Lemma 1.
 \hspace{.25in} $\Box$ 

\medskip

We now have as an immediate consequence the following.
\medskip
\begin{thm}{Suppose $M$ is an inner model of set theory and ${\Bbb R}^M$ is analytic.
Then either $\aleph _1^M$ is countable or all reals are in $M$.}
\end{thm}

\medskip

To extend Lemma 1 and consequently Theorem 3 to higher levels of 
the projective hierarchy we need the appropriate form of 
projective absoluteness in place of Shoenfield's theorem.
We first do this in the case of $\bf \Sigma$$_2^1$ sets.

\medskip

\begin{lemma}{Let $a$ be a real such that $a^{\#}$ exists and
assume that $A$ is a $\Sigma ^1_2(a)$ set such that
sup$\{ \omega _1^{CK,x}: x\in A\}=\aleph _1$. Then every real
is hyperarithmetic in a quadruple of elements of $A$.}
\end{lemma}

\noindent PROOF: Supposet $A$ is defined by a $\Sigma _2^1(a)$ formula
$\varphi (x,a)$. Following the proof of Lemma 1 we have to show that 
if $G$ is $V$-generic over $Coll(\aleph _0,\aleph _1)$ then in
$V[G]$ sup$\{ \omega _1^{CK,x}:\varphi (x,a)\ \mbox{holds}\}>\aleph _1^V$.
Let $\alpha <\aleph _1^V$ be indiscernible for $L[a]$. 
In $V$ pick an $L[a]$-generic filter $G_{\alpha}$ over $Coll(\aleph _0,\alpha)$.
This can be done since $\aleph _1^V$ is inaccessible in $L[a]$.
In $L[a,G_{\alpha}]$ pick a linear ordering  $R$ on $\omega$ such that
$(\omega,R)$ is isomorphic to $(\alpha, <)$. 
The formula which says that there exists $x$ such that 
$\varphi (x,a)$ holds and such that $\omega _1^{CK,x}>\alpha$ is
$\Sigma ^1_2(a,R)$ and is true in $V$. By Shoenfield's absoluteness
theorem it is true in $L[a,G_{\alpha}]$ as well.
Since $G_{\alpha}$ can be chosen to contain any condition in
$Coll(\aleph _0,\alpha)$ it follows that the maximal condition
in $Coll(\aleph _0,\alpha)$ forces the above statement.
Since both $\alpha$ and $\aleph _1^V$ are indiscernibles over $L[a]$
it follows that the maximal condition in $Coll(\aleph _0,\aleph _1^V)$
forces over $L[a]$  that there is $x$ such that $\varphi (x,a)$
holds and $\omega _1^{CK,x}>\aleph _1^V$. 

As in the proof of Lemma 1 we show that if $H$ is $V$-generic
over $Coll(\aleph _0,\aleph _2)$ then in $V[H]$ for any real 
$r$ there are reals  $x_1,x_2,x_3,y$ all satisfying $\varphi (x,a)$
and such that $r$ is $\Delta ^1_1(x_1,x_2,x_3,y)$. 
The existence of such quadruple is $\Sigma ^1_2(a,r)$
so if $r$ is in $V$ it follows, by Shoenfield's theorem again,
there such a quadruple exists already in $V$. 
\hspace{.25in} $\Box$ 
\medskip

\begin{thm}{Assume $x^{\#}$ exists, for every real $x$.
If $M$ is an inner model of set theory such that $\aleph _1^M$
is uncountable and ${\Bbb R}^M$ is $\bf \Sigma$$_2^1$ then all
reals are in $M$.  \hspace{.25in} $\Box$ }
\end{thm}

For an integer $n$ and an infinite cardinal $\kappa$ let us
say that a universe $V$ satisfies $\Sigma ^1_n$-absoluteness
for posets of size $<\kappa$ if whenever ${\cal P}$ is 
a forcing notion of size $<\kappa$ and 
in $V^{\cal P}$ ${\cal Q}$ is a forcing notion of size $< \kappa$
then for any $\Sigma ^1_n$ formula $\varphi$ with 
parameters from $V^{\cal P}$, $\varphi$ holds in $V^{{\cal P}\star{\cal Q}}$
if and only if it holds in $V^{\cal P}$. 
Woodin has shown that assuming the existence of $n$ Woodin cardinals
with a measurable cardinal above then $\Sigma _{n+3}^1$ absoluteness holds
for posets of size less than the first Woodin cardinal. 
The analogous proof to Lemma 1 goes through for $\bf \Sigma$$^1_{n+2}$
sets under this assumption.
Therefore we have the following. 

\begin{thm}{Assume the existence of $n$ Woodin cardinals 
with measurable above. If $M$ is an inner model of set theory 
such that $\aleph _1^M$ is uncountable and ${\Bbb R}^M$ 
is  a $\bf \Sigma$$^1_{n+2}$ set then all reals are in $M$.}
\end{thm}

\section{Adding perfect sets of ground model reals}

In this section we show that the conclusion of Theorem 3 cannot be strengthened to say 
that either all reals are in $M$ or the continuum of $M$ is countable. We also show
that it is possible to have an inner model of set theory $W$ such that $\aleph _1^W=\aleph _1$,
$W$ contains a perfect set of reals, and not all reals are in $W$.
We start with the following.

\begin{thm}{{\em (CH)} Suppose there is a club in $\omega _1$ consisting of
ordinals of uncountable cofinality in $L$. Then there is an $L$-generic
filter $G$ for adding $\omega _1^V$ many Cohen reals to $L$ such that 
the reals of $L[G]$ are an $F_{\sigma}$ set in $V$.}
\end{thm}

\medskip

\noindent PROOF: Let $C$ be a club in $\omega _1$ consisting of ordinals
of uncountable cofinality in $L$. Let $P$ be a perfect subset of 
$2^{\omega}$ such that any finite subest of $P$ consists of mutually generic
Cohen reals over $L$. Fix a recursive partition of $\omega$ into 
infinitely many disjoint infinite sets $\{ A_i:i<\omega \}$ and for each
$i<\omega$ fix a recursive partition $\{ A_{i,j}:j<\omega \}$ of $A_i$ 
into infinitely many disjoint infinite pieces. For each each $d\in 2^{\omega}$
let $d_i$ be the real obtained by restricting $d$ to $A_i$ and transfering
it to $2^{\omega}$ using the order preserving bijection between $A_i$
and $\omega$. Let $d_{i,j}$ be obtained by restricting $d$ 
to $A_{i,j}$ and transfering to $2^{\omega}$ in a similar fashion.

Construct the generic $G$ by constructing an $L$-generic 
filter $G_{\alpha}$ over ${\cal C}(\alpha)$ by induction on $\alpha \in C$.
 The requirements are that for each 
$\alpha \in C$ there exists a countable subset $S_{\alpha}$ of $P$ such that

\begin{enumerate}

\item for all $\beta < \alpha$ $G_{\alpha}(\beta)=d_{i,j}$, 
for some $d\in S_{\alpha}$, and some
$i,j<\omega$,

\item for all $d\in S_{\alpha}$ and for all $i,j<\omega$ there is $\beta <\alpha$
such that $G_{\alpha}(\beta)=d_{i,j}$,

\item the set of reals of $L[G_{\alpha}]$ is the union of the sets of 
reals in $L[s]$, where $s$ is a finite sequence of members of 
$\{d_i: d\in S_{\alpha}\ \mbox{and}\ i<\omega \}$.

\end{enumerate}

Since every $\alpha \in C$ has uncountable cofinality in $L$ genericity and
these conditions are preserved at a stage $\delta$ which is a limit point of $C$ by using
$S_{\delta}=\bigcup \{ S_{\alpha}:\alpha <\delta \}$. 
We now verify the successor step. 
Let $G_{\alpha}$ and $S_{\alpha}$ be given. By condition 3. any finite 
subset of $P\setminus S_{\alpha}$ consists of mutually generic Cohen reals
over $L[G_{\alpha}]$. Let $\alpha ^*$ be the next element of $C$ above $\alpha$.
Let $\{ X_i :i<\omega \}$ be an increasing sequence of subsets of $[\alpha ,\alpha ^*)$
such that each $X_i \in L[G_{\alpha}]$, $X_i$ is countable in $L[G_{\alpha}]$,
and such that if $Y\subset [\alpha ,\alpha ^*)$ is countable in $L[G_{\alpha}]$
then $Y\subseteq X_i$, for some $i<\omega$. Moreover arrange that 
$X_{i+1}\setminus X_i$ is infinite, for each $i$. 
Fix any $d\in P\setminus S_{\alpha}$. It is routine to construct 
$G^*$ satisfying 1. and 2. for $S_{\alpha ^*}=S_{\alpha} \cup \{ d \}$
 and such that for all $i$
$$
L[G_{\alpha}][g_i]=L[G_{\alpha}[d_i]
$$
where $g_i=G_{\alpha ^*} \restriction (X_i\setminus X_{i-1})$.
Then condition 3. follows.  

Assuming CH we can easily arrange that 
$P=\bigcup \{ S_{\alpha}: \alpha <\omega _1\}$. Thus the set of reals
in $L[G]$ is exactly the union of the reals of $L[s]$, where 
$s$ is a finite sequence of elements of 
$\{ d_i: d\in P\ \mbox{and}\ i<\omega \}$. 
Since there are only countably many terms for reals in Cohen extensions
which are in $L$ and $P$ is compact, it follows that this set
is $F_{\sigma}$.
\hspace{.25in} $\Box$

\medskip

To obtain a model satisfying the assumptions of Theorem 6 we can start with
a model of $V=L$, collapse $\aleph _1$ to $\aleph _0$ and then shoot
a club through the set of ordinals $<\aleph _2^L$ of uncountable cofinality in $L$.
Thus we have the following.

\medskip 

\begin{thm}{There is a pair $V$ and $W$ of generic extensions of $L$ such 
that $W\subseteq V$, the reals of $W$ form an uncountable $F_{\sigma}$ 
set in $V$, and $V$ and $W$ do not have the same reals. \hspace{.25in} $\Box$}
\end{thm}

\medskip

The following result says that we can have an inner model of set theory for 
which Prikry's question has a positive answer. 

\medskip

\begin{thm}{ Assume ZFC. Then there is a pair $(W,V)$ of generic 
extensions of $L$ such that $W\subseteq V$, $\aleph _1^W=\aleph _1^V$,
and $V$ contains a perfect $P$ set of $W$-reals which is not in $W$.}
\end{thm}

We will need the following lemma (cf. Theorem 1 from [SW]).

\medskip 

\begin{lemma}{There is a generic extension $V_0$ of $L$ such that
 $\aleph _1^{V_0}=\aleph _1^L$, and $V_0$ contains a club $C$ in
$\aleph _3^L$ consisting of ordinals of uncountable cofinality in $L$.}
\end{lemma}

\medskip

\noindent PROOF: $V_0$ will be obtained as a two step forcing extension
of $L$. Let ${\cal N}$ be the following version of Namba forcing.
Conditions in ${\cal N}$ are subtrees $T$ of $\omega _2^{<\omega}$ 
such that for every $s\in T$ the set $\{ t\in T: s\subseteq t\}$ 
has cardinality $\aleph _2$. The partial ordering is defined in
the natural way: $R\leq T$ if and only if $R\subseteq T$. 
For a node $s\in T$ we let $T_s=\{ t\in T: t\subseteq s\ \mbox{or}\ s\subseteq t\}$.
Then ${\cal N}$ preserves $\aleph _1$, changes the cofinality of $\aleph _2$
to $\aleph _0$, and  collapses the cardinality of $\aleph _3$ to $\aleph _1$.
Define in $L$ the set $S=\{ \alpha <\omega _3: cof(\alpha)=\omega _2\}$.
Suppose now that $G$  is $L$-generic over ${\cal N}$.

\medskip

\noindent {CLAIM:} $S$ remains stationary in $L[G]$.

\medskip

\noindent PROOF: Working in $L$ let a name $\dot{C}$ for a club in $\omega _3$
and a condition $T\in {\cal N}$ be given. 
Fix a sufficiently large regular 
cardinal $\theta$ and take an elementary submodel $M$ of $H_{\theta}$ of cardinality
$\aleph _2$ containing $\dot{C}$ and $T$ such that $M\cap \omega _3
=\delta \in S$.
 
By shrinking if necessary we may assume that every node in $T$ has either $1$ 
or $\aleph _2$ immediate extensions. Fix a strictly increasing sequence 
$\langle \delta _{\xi}:\xi <\omega _2\rangle$ of ordinals
 converging to $\delta$.
We  build by a fusion argument a condition $R\leq T$ such that 
$R\forces \delta \in \dot{C}$.
Set $R_0=T$. Let $s$ be the stem of $T$.
For each $\xi <\omega _2$ such that $s\ \wh  \xi \in T$ 
the condition $T_{s\ \wh  \ \xi}$ belongs to ${\cal N}\cap M$.
By elementarity and the fact that $\dot{C}$ is forced to be unbounded in $\omega _3$
there is a condition $Q_{s,\xi}\leq T_{s\ \wh  \xi}$ 
such that $Q_{s,\xi}\in {\cal N}\cap M$ 
and for some $\delta _{\xi} <\gamma <\delta$ 
$Q_{s,\xi}\forces \gamma \in \dot{C}$.
Let 
$$
R_1=\bigcup \{ Q_{x,\xi}:\xi <\omega _2\ \mbox{and}\ {s\  \wh \xi} \in T\}.
$$
Now given $R_n$ let $L_n$ be the set of nodes of $R_n$
 which are $\aleph _2$-splitting
and have exactly $n$ $\aleph _2$-splitting nodes below them.
For each $t\in L_n$ we have $(R_n)_t\in M$ so, by a similar argument,
 for each $\xi <\omega _2$ such that $t\ \wh  \xi \in R_n$ 
we can pick $Q_{t,\xi} \leq R_{t\ \wh \xi}$
with $Q_{t,\xi}\in M$ such that for some $\delta _{\xi}<\gamma <\delta$ 
$Q_{t,\xi}\forces \gamma \in \dot{C}$. 
Then we let 
$$
R_{n+1}=\bigcup \{ Q_{t,\xi}:t \in L_n\ \mbox{and}\ t\ \wh \xi \in R_n\}.
$$
Finally let $R=\bigcap \{ R_n:n<\omega \}$. 
Then $R\in {\cal N}$ and if $t$ is an $\aleph _2$-splitting node of $R$ 
it follows that for every $\xi <\omega _2$ such that $t\ \wh \xi \in R$ 
there is $\delta _{\xi} <\gamma < \delta$ such that 
$R_{t\ \wh  \xi}\forces \gamma \in \dot{C}$. This implies that 
$R\forces \delta \in \dot{C}\cap S$, as required.
\hspace{.25in} $\Box$

\medskip

Now if $G$ is $L$-generic over ${\cal N}$ let in $L[G]$ ${\cal Q}$ be 
the standard poset for shooting a club through $S$ with countable conditions. 
Then if $C$ is the generic club it consists of ordinals of $L$-cofinality
$\aleph _2$. 
\hspace{.25in} $\Box$

\medskip

\noindent PROOF of Theorem 7:  
For any index set $I$ let ${\cal C}(I)$ denote the standard poset for
adding $I$ Cohen reals.
Let ${\cal P}$ be the poset for adding a perfect set of mutually generic
Cohen reals, that is a perfect set $P_g$ of reals such that for
any 1-1 sequence $\bar{b}$ of length $n$ of  members of $P_g$
$\bar{b}$ is $V$-generic for ${\cal C}(n)$. 
A condition $\sigma$ belongs to ${\cal P}$ if there is
an integer $m =m(\sigma)$
such that $\sigma$ is an initial segment of $\{ 0,1\}^{\leq m}$ 
with the property that every $s\in \sigma$ has an extension in $\sigma$ 
of height $m$. Say that $\tau \leq \sigma$ iff 
$\tau \restriction \{ 0,1\}^{\leq m(\sigma)} =\sigma$.
Thus, in terms of forcing, ${\cal P}$ is equivalent to the standard poset for 
adding a single Cohen reals. 
If $g$ is $V$-generic for ${\cal P}$ then $T_g=\bigcup g$
is a perfect tree. Let $P_g=[T_g]$ denote the set of all infinite
branches of $T_g$ as computed in the model $V[g]$.

Let now $V_0$ be the generic extension of $L$ as in Lemma 2. 
We shall force over $V_0$ with the poset ${\cal C}(\omega _3^L)\times {\cal P}$.
Note that this poset is equivalent to ${\cal C}(\omega _1)$.
Suppose $G\times g$ is $V_0$-generic for ${\cal C}(\omega _3^L)\times {\cal P}$.
Then we can identify $G$ with an $\omega _3^L$-sequence 
$\langle G({\xi}):\xi <\omega _3^L\rangle$ of Cohen reals.
Let $P=P_g^{V_0[G\times g]}$ denote $[T_g]$ as computed in
the model $V_0[G\times g]$. Note that since the forcing notion
${\cal P}$ is the same whether defined in $V_0$ or $V_0[G]$ 
we conclude that the reals in $P$ are mutually Cohen generic
over $V_0[G]$. 

In $V_0$ fix a club $C$ in $\omega _3^L$ consisting of ordinals 
of uncountable cofinality in $L$. Note that for any $X\in L$ which
is countable in $L$ $X\cap C$ is finite. 
In $V_0$ fix an enumeration $\{ r_{\alpha}:\alpha <\omega _1\}$
of $P$ and an increasing enumeration $\{ \gamma _{\alpha}:\alpha <\omega _1\}$
of $C$. 
We now define an $\omega _3^L$-sequence of reals $G^*$ as follows. 
If $\gamma =\gamma _{\alpha}$ for some $\alpha <\omega _1$ then 
let $G^*(\gamma)=r_{\alpha}$, otherwise let $G^*(\gamma)=G(\gamma)$.

\medskip 

\noindent CLAIM: $G^*$ is $L$-generic over ${\cal C}(\omega _3^L)$.

\medskip

\noindent PROOF: 
Since ${\cal C}(I)$ has the ccc for any $I$ it suffices to show that
for any $I\subseteq \omega _3^L$ which is countable in $L$ 
$G^*\restriction I$ is $L$-generic over ${\cal C}(I)$.
Fix such $I$. By the property of the club $C$ it follows that 
$I\cap C$ is finite. Let $F\subseteq \omega _1$ be finite such 
that $I\cap C\subseteq \{ \gamma _{\alpha}:\alpha \in F\}$. 
Now $G^*\restriction (I\setminus F)=G\restriction (I\setminus F)$
and the sequence $\langle r_{\alpha}:\alpha \in F\rangle$ is 
$L[G]$-generic over ${\cal C}(F)$. It follows that $G^*\restriction I$
is $L$-generic over ${\cal C}(I)$. 
\hspace{.25in} $\Box$

\medskip
 
Now let $W=L[G^*]$ and  $V=V_0[G\times g]$. By the definition of
$G^*$ we have that $P_g\subseteq W$. We claim that $T_g$ does not
belong to $W$. Otherwise there would be a countable $I\subseteq \omega _3^L$
such that $I\in L$ and $T_g\in L[G^*\restriction I]$. Since $T_g$ is a perfect
tree it would have infinitely many branches in $L[G^*\restriction I]$.
Since $I\cap C$ is finite there would exist $\alpha \in \omega _1$
such that $\gamma _{\alpha}\notin I$ and $r_{\alpha}\in L[G^*\restriction I]$.
This contradicts the fact that $r_{\alpha}$ is Cohen generic over $L[G^*\restriction I]$.
\hspace{.25in} $\Box$

\medskip

\section{Submodels of $L({\Bbb R})$ under AD}

We now show how the coding techniques introduced in previous sections can be applied 
in the context of AD + $V=L({\Bbb R})$. The following result implies that under 
this assumption the property of being a cardinal below $\theta$ is $\Delta _1$
over $L({\Bbb R})$. 

\medskip

\begin{thm}{Assume {\em AD + }$V=L({\Bbb R})$. If $M$ is an inner model of {\em ZF}
containing a Souslin prewellordering of reals of length $\aleph _1^V$ then all reals
are in $M$.}
\end{thm}

Some assumptions on the prewellordering in Theorem 9 are necessary. We show the following.

\begin{thm}{{\em (ZF)} Assume there is a nonconstructible real. Then there is a transitive inner model
$M$ of {\em ZF} in which there is a prewellordering of the reals of length $\omega _1^V$ 
and such that not all reals belong to $M$.}
\end{thm}

\medskip

\noindent PROOF: Let $\langle c_i: i<\omega \rangle$ be a sequence of mutually generic Cohen 
reals over $L$. Let $S$ be the set of reals constructible from  finitely many of the $c_i$'s
and let $T$ be the set of Turing degrees of the $c_i$'s. Then $L(S,T)$ is a symmetric
extension of $L$. For an ordinal $\delta$  in  $L(S,T)$ consider the partial ordering
 ${\cal Q}$ for adding a map from $T$ to $\delta$ with finite conditions. 
Thus members of ${\cal Q}$ are finite partial functions from $T$ to $\delta$ and 
the ordering is reverse inclusion. We can identify the generic filter $G$ with a prewellordering
$\leq _G$ of $T$ where $\tau \leq _G \sigma$ iff $(\bigcup G)(\tau)\leq (\bigcup G)(\sigma)$.

\medskip

\noindent CLAIM: If $\leq$ is {\em any} prewellordering of $T$ of length $\delta$ for which
the induced equivalence classes are infinite then $\leq$ is $L(S,T)$-generic over ${\cal Q}$.
Moreover $L(S,T)[\leq]$ and $L(S,T)$ have the same reals.

\medskip

\noindent PROOF: Let $\leq$ be any prewellordering of $T$ satisfying the requirements of the claim
and let $H$ be the corresponding filter in ${\cal Q}$. 
Then $h=(\bigcup H):T\rightarrow \delta$ and  $h^{-1}(\xi)$ is infinite, for every $\xi <\delta$.
Let $D\in L(S,T)$  be a dense subset of ${\cal Q}$. We have to show that $D\cap H \neq \emptyset$.
There is $n<\omega$ such that $D$ is definable in $L(S,T)$ from parameters 
$\{ c_1,\ldots,c_n\} \cup \{ S,T\}$. For each $i$ let $d_i$ be the Turing degree 
of $c_i$. 
Let $F=\{ d_1,\ldots,d_n\}$ and let $p=h\restriction F$. 
Then $p\in {\cal Q}$. Using the density of $D$ find $q\leq p$ such that $q\in D$.
We may assume without loss of generality that for some $m\geq n$ 
dom$(q)=\{ d_1,\ldots,d_m\}$. By the property of $h$ we can find a 1-1 function 
$f:m\setminus n \rightarrow \omega \setminus n$ such that for all $i\in [n,m)$ $q(d_i)=h(d_{f(i)})$.
Let $q^*=h\restriction (F \cup \{ d_{f(j)}: n\leq j <m \})$.
 We show that $q^*\in D$.
To see this fix a recursive permutation $\varphi$ of $\omega$ 
extending $(id\restriction n) \cup f$. $\varphi$ induces a permutation of $\{ c_i:i<\omega \}$
which in turn induces an automorphism $\varphi ^*$ of $L(S,T)$ which fixes
$c_1,\ldots,c_n$, and each Turing degree in $T$. Then $\varphi (D)=D$ and $\varphi ^*(q)=q^*$. 
From this it follows that $q^* \in D$, as required.

To prove that ${\cal Q}$ does not add any reals to $L(S,T)$ let $H$ and $h$ 
be as above and suppose $\dot{r}$
is a ${\cal Q}$-name for a real. Then as before there is $n$ such that $\dot{r}$
is definable from $\{ c_1,\ldots,c_n\} \cup \{ S,T\}$. Let $F=\{ d_1,\ldots ,d_n\}$
and $p=h\restriction F$. 
Let $m<\omega$ and suppose a condition $q\leq p$ decides the value of $\dot{r}(m)$.
Then as in the previous argument there is a condition $q^*\in H$ such that 
some automorphism of $L(S,T)$ fixes $\dot{r}$ and maps $q$ to $q^*$.
Thus $q^*$ forces the same information about $\dot{r}(m)$ as $q$. 
This implies that $p$ forces that $\dot{r}$ is in $L(S,T)$ as desired.
\hspace{.25in} $\Box$

\medskip

 To finish the proof of Theorem 10 notice that we may assume that $\omega _1^L$
is countable since otherwise we can take $M=L$. Let $P$ be a perfect set of 
mutually generic Cohen reals over $L$. Let $S$ be the set of reals constructible from
finitely many members of $P$ and let $T$ be the set of Turing degrees of the $c_i$'s.
Let $\leq$ be any prewellordering of $T$ of length $\omega _1^V$ whose induced 
equivalence classes are infinite. Then $\leq$ will be generic over $L(S,T)$.
To see this go to a generic extension of the universe in which the continuum of $V$
is countable and apply the claim.  Let $M=L(S,T)[\leq]$. Then by applying the second
part of the claim $M$ and $L(S,T)$ have the same reals and therefore not all 
reals are in $M$.  Therefore $M$ satisfies the conclusions of the theorem. 
\hspace{.25in} $\Box$

\vspace{.25in}

Carnegie Mellon University, Pittsburgh PA 15213 and Universit\'e Paris VII,
Jussieu, 75251 Paris

\medskip

University of California, Berkeley, CA 94720

\end{document}